\documentclass[11pt, oneside]{article}   	
\usepackage{amsmath, amsthm, amsfonts, amssymb}
\usepackage{graphicx}
\usepackage[colorlinks,citecolor=blue,urlcolor=blue]{hyperref}
\newcommand\Ex{{\mathbb E}}

\newcommand\X{{\mathbb X}}

\newcommand\cC{{\mathcal C}}

\newcommand\cP{{\mathcal P}}

\newcommand\N{{\mathbb N}}

\newcommand\R{{\mathbb R}}
\newcommand\C{{\mathbb C}}

\newcommand\cc{{c}}

\newcommand\bra[1]{\langle #1 \rangle}

\DeclareMathOperator{\Beta}{Beta}
\DeclareMathOperator{\supp}{supp}
\DeclareMathOperator{\tr}{trace}

\newtheorem{theorem}{Theorem}[section]

\newtheorem{lemma}[theorem]{Lemma}

\theoremstyle{definition}

\theoremstyle{remark}

\title{Beta Jacobi ensembles and associated Jacobi polynomials}

\author{Hoang Dung Trinh\footnote{Faculty of Mathematics Mechanics Informatics, University of Science, Vietnam National University, Hanoi, Vietnam.
\newline Email: thdung.hus@gmail.com} 
\and
Khanh Duy Trinh \footnote{Global Center for Science and Engineering, Waseda University, Japan.
\newline
Email: trinh@aoni.waseda.jp 
} 
}

\begin{document}
\maketitle

\begin{abstract}
Beta ensembles on the real line with three classical weights (Gaussian, Laguerre and Jacobi) are now realized as the eigenvalues of certain tridiagonal random matrices. The paper deals with beta Jacobi ensembles, the type with the Jacobi weight. Making use of the random matrix model, we show that in the regime where $\beta N \to const \in [0, \infty)$, with $N$ the system size, the empirical distribution of the eigenvalues converges weakly to a limiting measure which belongs to a new class of probability measures of associated Jacobi polynomials. This is analogous to the existing results for the other two classical weights. We also study the limiting behavior of the empirical measure process of beta Jacobi processes in the same regime and obtain a dynamical version of the above.

\medskip

		
\medskip
	
\end{abstract}

\section{Introduction}
Real beta ensembles are ensembles of real particles distributed according to the following joint probability density function
\begin{equation}\label{joint-density}
	Z\cdot \prod_{i < j}|\lambda_j - \lambda_i|^\beta \prod_{l = 1}^N w(\lambda_l) ,\quad (\beta > 0),
\end{equation}
where $w(\lambda)\ge 0$ is a weight function and $Z$ is the normalizing constant. They play important roles in many areas such as mathematical physics, statistical mechanics, random matrix theory, multivariate statistical theory and representation theory.
With three classical weights 
\begin{equation}\label{weights}
	w(\lambda) = 
	\begin{cases}
		e^{-\lambda^2/2},\quad &\lambda \in (-\infty, \infty), \quad \text{Gaussian,}\\
		\lambda^\alpha e^{-\lambda}, \quad &\lambda \in (0, \infty),\quad (\alpha > -1),\quad \text{Laguerre,}\\
		\lambda^a (1 - \lambda)^b, \quad &\lambda \in (0,1), \quad(a>-1,b>-1),\quad \text{Jacobi},
	\end{cases}
\end{equation}
the ensembles are called Gaussian beta ensembles, beta Laguerre ensembles and beta Jacobi ensembles, respectively. The three classical beta ensembles with specific values of $\beta$ were originally known as the eigenvalues of random matrices: Gaussian orthogonal/unitary/symplectic ensembles (Gaussian weight, $\beta = 1,2,4$), Wishart/Laguerre matrices (Laguerre weight, $\beta = 1,2$), and `double Wishart' matrices (Jacobi weight, $\beta = 1$). However, they are now realized as the eigenvalues of certain tridiagonal random matrices for any $\beta > 0$ \cite{DE02, Killip-Nenciu-2004}.

Study the limiting behavior of the empirical distribution of the eigenvalues is a very first problem in random matrix theory. Under some mild conditions on the weight $w$, the empirical distribution converges weakly to a limiting probability measure $\mu_w$ (called the equilibrium of the system), almost surely. The convergence to the equilibrium $\mu_w$ holds even when the parameter $\beta$ varies but satisfies $\beta N \to \infty$. What happens when $\beta N $ stays bounded? It turns out that the empirical distribution converges to a different limit $\mu_c$ in the regime where $\beta N \to 2c \in [0, \infty)$ \cite{G-Zelada-2019, Lambert-2019, Nakano-Trinh-2020}. Note that all those results have been proved by analyzing the joint density of the beta ensembles.

For Gaussian beta ensembles and beta Laguerre ensembles, based on their tridiagonal random matrix model and a duality relation between $\beta$ and $4/\beta$, the limiting measure $\mu_c$ in the regime $\beta N \to 2c$ can be calculated explicitly. It was shown in \cite{DS15, Trinh-Trinh-2019} that in the Gaussian case (resp.\ Laguerre case)  the limit $\mu_c$ belongs to a family of probability measures of associated  Hermite polynomials (resp.\ associated Laguerre polynomials). Here Hermite polynomials (resp.\ Laguerre polynomials) are orthogonal polynomials with respect to the Gaussian weight (resp. the Laguerre weight). Their associated orthogonal polynomials obtained by shifting the coefficients in the three term recurrence relation were studied in 1980s in \cite{Askey-Wimp-1984, Ismail-et-al-1988}. These motivate this study to see if an analogous phenomenon happens for beta Jacobi ensembles. We find out that the limiting measure $\mu_c$ in the Jacobi case belongs to a new family of associated Jacobi polynomials which is slightly different from the two existing ones \cite{Ismail-Masson-1991}.

Let us introduce some preliminaries before stating our main result. Let $\mu$ be a nontrivial ($\mu$ is not supported on finite points) probability measure on $\R$ with all finite moments. Then the set $\{1, x, x^2, \dots\}$ is linearly independent in $L^2(\R, \mu)$. Let $P_n = x^n + \text{lower orders}, P_0=1$ be the orthogonal polynomials resulting from the Gram Schmidt orthogonalization process applying to that set. Then the polynomials $\{P_n\}$ satisfy a three term recurrence relation
\begin{equation}\label{three-term}
	xP_n = P_{n+1} + a_{n+1} P_n + b_n^2 P_{n-1}, \quad n \ge 0, 
\end{equation}
for coefficients $a_n \in \R, b_n >0, n\ge 1, (b_0:=0)$. 
Moreover, $\{p_n = (b_1 \cdots b_n)^{-1}P_n\}_{n \ge 0}$ become an orthonormal system and the linear transformation of multiplication by $x$  in $L^2(\R, \mu)$ has the following matrix 
\begin{equation}\label{Jacobi-matrix}
	J=\begin{pmatrix}
		a_1 	&b_1\\
		b_1		&a_2	&b_2\\
		&\ddots&\ddots&\ddots
	\end{pmatrix},
\end{equation}
namely, $J \vec{p} = x\vec{p}$, where $\vec{p} =(p_0, p_1, \dots)^t$. The matrix $J$ is called the Jacobi matrix of the probability measure $\mu$.

The inverse problem is to find a probability measure which orthogonalizes polynomials $\{P_n\}$ satisfying a three term recurrence relation \eqref{three-term} for given  $\{a_n \in \R, b_n>0\}_{n \ge 1}$, or equivalently, for given infinite Jacobi matrix $J$. This is a classical problem related to the Hamburger moment problem. Any probability measure $\mu$ satisfying  
\begin{equation}\label{definition-of-spectral-measure}
	\int x^k d\mu(x) = J^k(1,1), \quad k = 0,1,\dots,
\end{equation}
is a solution. As a measure $\mu$ satisfying the moments condition \eqref{definition-of-spectral-measure} always exists, it is, in general, not unique. The uniqueness of the moment problem which is equivalent to the essential self-adjointness of the operator $J$ on $\ell^2(\N)$ holds under a useful sufficient condition $\sum b_n^{-1} = \infty$ \cite[Corollary 3.8.9]{Simon-book-2011}. The measure $\mu$ is called the spectral measure of $J$, or the probability measure of the polynomials $\{P_n\}$ in case of unicity.

Jacobi polynomials are orthogonal polynomials with respect to the Jacobi weight $x^a(1-x)^b dx, 0 < x < 1$, or the probability measure $const \times x^a(1-x)^b dx, 0 < x < 1$. From the three term recurrence relation of Jacobi polynomials, we deduce that the Jacobi matrix of Jacobi polynomials is given by
\begin{align*}
	J_{Ja} &= 	\begin{pmatrix}
		\sqrt{\lambda_0}	\\
		\sqrt{\mu_1}	&\sqrt{\lambda_1}	\\
		&\ddots	&\ddots	
	\end{pmatrix}
	\begin{pmatrix}
		\sqrt{\lambda_0}	&\sqrt{\mu_1}	\\
			&\sqrt{\lambda_1}	&\sqrt{\mu_2}	\\
			&&\ddots	&\ddots	
	\end{pmatrix},\\
	&\qquad
	\begin{cases}
	\lambda_n = \frac{n + a+ 1}{2n + a + b + 2}\frac{n + a + b + 1}{2n+ a + b + 1}, &n \ge 0,\\
	\mu_n =  \frac{n}{2n + a + b + 1}\frac{n  + b}{2n + a + b}, & n \ge 1.
	\end{cases}
\end{align*}

Our main result is as follows.
\begin{theorem}\label{thm:main-intro}
	Let $L_N = N^{-1} \sum_{i = 1}^N \delta_{\lambda_i}$ be the empirical distribution of the beta Jacobi ensemble 
\[
Z\cdot \prod_{i < j}|\lambda_j - \lambda_i|^\beta \prod_{l = 1}^N \lambda_l^a (1-\lambda_l)^b, \quad \lambda_i \in (0,1).
\]
Here $\delta_\lambda$ denotes the Dirac measure at $\lambda$. Let $a, b > -1$ be fixed. Then in the regime where $\beta N \to 2c \in [0, \infty)$, the empirical distribution $L_N$ converges weakly to a limiting measure $\nu_{a, b, c}$, almost surely. Here $\nu_{a,b,c}$ is the spectral measure of the following Jacobi matrix
\begin{align*}
		J_c &= 	\begin{pmatrix}
		\sqrt{\hat\lambda_0(c)}	\\
		\sqrt{\mu_1(c)}	&\sqrt{\lambda_1(c)}	\\
		&\ddots	&\ddots	
	\end{pmatrix}
	\begin{pmatrix}
		\sqrt{\hat \lambda_0(c)}	&\sqrt{\mu_1(c)}	\\
			&\sqrt{\lambda_1(c)}	&\sqrt{\mu_2(c)}	\\
			&&\ddots	&\ddots	
	\end{pmatrix},\\
	&\qquad
	\begin{cases}
	\hat\lambda_0(c) = \frac{c+a+1}{2c+a+b+2},\\
	\lambda_n(c) = \frac{n + c + a+ 1}{2n + 2c + a + b + 2}\frac{n + c + a + b + 1}{2n + 2c + a + b + 1}, \\
	\mu_n(c) = \frac{n + c}{2n + 2c + a + b + 1}\frac{n + c + b}{2n + 2c + a + b}, & n \ge 1.
	\end{cases}
\end{align*}
\end{theorem}

To show this result, we will make use of the tridiagonal random matrix model together with some ideas already used in the study of Gaussian beta ensembles and Laguerre beta ensembles. 
Note that when $c = 0$, the Jacobi matrix $J_c$ coincides with $J_{Ja}$. For $c = 1, 2, \dots,$ except the first $\hat \lambda_0$, the parameters in $J_c$ are shifted from those in $J_{Ja}$, that is, $\lambda_n(c) = \lambda_{n + c}, \mu_n(c) = \mu_{n + c}$. Thus, we called the model with matrix $J_c$ associated Jacobi polynomials. This new model is slightly different from the two existing ones of associated Jacobi polynomials. Since arguments are almost the same as those used in the Gaussian case and the Laguerre case, we do not give detailed proofs. The main contribution of this paper is to calculate the explicit formula for the new model of associated polynomials, which will be given in Sect.\ 4.

The paper also establishes a dynamical version of that result.  
The paper is organized as follows. In Sect.\ 2, we introduce the tridiagonal random matrix model for beta Jacobi ensembles and then prove Theorem~\ref{thm:main-intro}. Section 3  shows the result on the dynamical version. Finally, Sect.\ 4 deals with  associated Jacobi polynomials.

\section{Beta Jacobi ensembles}

A tridiagonal random matrix model for beta Jacobi ensembles was introduced in \cite{Killip-Nenciu-2004} as follows. Denote by $\Beta(a,b)$ the beta distribution with parameter $a, b > 0$.
Let $p_1, \dots, p_N$ and $q_1, \dots, q_{N - 1}$ be independent random variables with 
\begin{align*}
	p_n &\sim \Beta((N - n) \kappa + a + 1, (N - n) \kappa + b + 1),\\
	q_n &\sim \Beta((N - n) \kappa, (N - n - 1) \kappa + a + b + 2),
\end{align*}
where $a, b > -1$ and $\kappa = \beta/2$. 
Let 
\begin{align*}
	s_n &= \sqrt{p_n(1 - q_{n - 1})},\quad n = 1, \dots, N, \quad (q_0 = 0),\\
	t_n &= \sqrt{q_n(1 - p_n)}, \quad n = 1, \dots, N - 1.
\end{align*}
Then a tridiagonal random matrix
\[
	J_{N, \beta}(a, b) =  \begin{pmatrix}
		s_1	\\
		t_1	&s_2		\\
		&\ddots	&\ddots \\
		&& t_{N - 1}	& s_N
	\end{pmatrix}
	 \begin{pmatrix}
		s_1 	&t_1	\\
			&s_2		&t_2		\\
		&&\ddots	&\ddots \\
		&&& s_N
	\end{pmatrix}
\]
has the eigenvalues $(\lambda_1, \dots, \lambda_N)$ distributed as the beta Jacobi ensemble
\[
	Z \times \prod_{i < j} |\lambda_j - \lambda_i|^\beta \prod_{l = 1}^N \lambda_l^a (1 - \lambda_l)^b, \quad \lambda_i \in [0,1],
\]
with $Z$ being the normalizing constant. Note that for $\beta = 1$, the ensembles are the eigenvalues of double Wishart matrices \cite[Theorem 3.3.4]{Muirhead-book-1982}. Based on this random matrix model, the limiting behavior of the empirical distribution, for fixed $\beta$, was studied in \cite{Dumitriu-Paquette-2012}.

Denote by $L_N$ the empirical distribution of the eigenvalues $\{\lambda_i\}$,
\[
	L_N = \frac 1N \sum_{i = 1}^N \delta_{\lambda_i}.
\]
Let $m_k(N, \kappa, a, b)$ be the expected value of the $k$th moment of the empirical distribution $L_N$, that is, 
\[
	m_k(N, \kappa, a, b) = \Ex[\bra{L_N, x^k}] = \Ex\left[\frac1N \tr(J_{N, \beta}(a, b)^k)\right].
\]
Here we denote $\bra{\mu, f} = \int f d\mu$ for an integrable function $f$ with respect to a measure $\mu$.
It turns out that $m_k(N, \kappa, a, b)$ can also be expressed as 
\[
	m_k(N, \kappa, a, b) = \Ex [J_{N, \beta}(a, b)^k(1,1)].
\]
(This property holds for Gaussian beta ensembles, beta Laguerre ensembles as well (see \cite{Duy-2018}, for example).) Then from formulas for moments of the beta distribution, we see that $m_p$ can be defined for any $N, \kappa, a$ and $b$. Our next arguments are based on the following duality relation which is a consequence of a result in \cite[Appendix A]{Dumitriu-Paquette-2012} (see also \cite[Eq.\ (4.15)]{Forrester-2017}).

\begin{lemma}
The following relation holds
	\[
		m_k(N, \kappa, a, b) =  m_k(-\kappa N, 1/\kappa, - a/\kappa, -b/\kappa).
	\]
\end{lemma}

The duality relation suggests that the regime where $\kappa N \to c$ with fixed $a, b$ can be studied by considering the regime where $N$ is fixed, $\kappa \to \infty,$ and $a=A \kappa, b = B \kappa$, for fixed $ A, B > 0$. In a new regime, since
\begin{align*}
	p_{n} \sim \Beta((N - n)\kappa + A \kappa + 1, (N-n)\kappa + B \kappa + 1), \quad n = 1,2,\dots, N,\\
	q_{n} \sim \Beta((N - n)\kappa, (N - n - 1)\kappa + A \kappa + B \kappa + 2), \quad n = 1, 2, \dots, N - 1,
\end{align*}
it follows immediately from the limiting behavior of the beta distribution that as $\kappa \to \infty$,
\begin{align*}
	p_{n} &\to \frac{n - N - A}{2n - 2N - A - B}, \quad n = 1,2,\dots, N,\\
	q_{n} &\to \frac{n - N}{2n - 2N - A - B + 1}, \quad n = 1, 2, \dots, N - 1,
\end{align*}
and thus
\begin{align}
	s_1^2 &= p_1 \to \frac{1 - N - A}{2 - 2N - A - B},\label{S1}\\
	s_n^2 &=p_n (1 - q_{n - 1}) \to  \frac{n - N - A}{2n - 2N - A - B} \frac{n - N - A - B}{2n - 2N - A - B - 1}, \quad n \ge 2, \label{Sn}\\
	t_n^2 &= q_n(1 - p_n) \to  \frac{n - N}{2n - 2N - A - B + 1} \frac{n - N - B}{2n - 2N - A - B}, \quad n \ge 1.\label{Tn}
\end{align}

Next, we define 
\begin{align*}
	\hat\lambda_0(c) & =  \frac{c + a + 1}{2c + a + b + 2},\\
	\lambda_n(c) & =  \frac{n + c + a + 1}{2n + 2c + a + b + 2} \frac{n + c + a + b + 1}{2n + 2c + a + b + 1}, \quad n \ge 1,\\
	\mu_n(c) & =  \frac{n + c}{2n + 2c + a + b + 1} \frac{n + c + b}{2n + 2c + a + b}, \quad n \ge 1,
\end{align*}
by exchanging $N \leftrightarrow -\cc, A \leftrightarrow -a, B \leftrightarrow -b$ in the limits~\eqref{S1}, \eqref{Sn} and \eqref{Tn} and then form an infinite Jacobi matrix $J_{III}$
\begin{align*}
	J_{III}&=\begin{pmatrix}
		\sqrt{\hat\lambda_0}	\\
		\sqrt{\mu_1}	&\sqrt{\lambda_1}		\\
		&\ddots	&\ddots 
	\end{pmatrix}
	\begin{pmatrix}
		\sqrt{\hat\lambda_0} &	\sqrt{\mu_1}	\\
		&\sqrt{\lambda_1}	&\sqrt{\mu_2}		\\
		&&\ddots	&\ddots 
	\end{pmatrix}\\
	& =\begin{pmatrix}
		\hat\lambda_0 	& \sqrt{\hat\lambda_0 \mu_1}	\\
		 \sqrt{\hat\lambda_0 \mu_1}	&\lambda_1 + \mu_1		& \sqrt{\lambda_1 \mu_2}	\\
		 &\ddots	&\ddots	&\ddots	
	\end{pmatrix}.
\end{align*}
Here for simplicity, we have removed the dependence on $c$ in formulas.
Then the duality relation implies the following result. 
\begin{lemma}
	For fixed $a, b > -1$, as $\kappa N \to c \in [0, \infty)$,
	\[
		m_k(N, \kappa, a, b) \to (J_{III})^k(1,1).
	\] 
\end{lemma}

The convergence of the expected values, together with the tridiagonal random matrix model, implies the almost sure convergence of moments of the empirical distribution $L_N$. The reason is that  
\[
	\bra{L_N, x^k} = \frac1N \tr(J_{N, \beta}(a, b)^k) = \frac1N \sum_{i = 1}^N (J_{N, \beta}(a, b)^k)(i, i)
\]
is a sum whose summands $(J_{N, \beta}(a, b)^k)(i, i)$ and $(J_{N, \beta}(a, b)^k)(j, j)$ are independent, if $|i - j|$ is large enough. A detail is omitted because it is similar to arguments used in case of Gaussian beta ensembles \cite{Trinh-2019} and of beta Laguerre ensembles \cite{Trinh-Trinh-2019}.

Let $\nu_{a,b,c}$ be the spectral measure of the Jacobi matrix $J_{III}$, that is, a unique measure satisfying 
\[
	\int x^k d\nu_{a,b,c} = (J_{III})^k(1,1), \quad k = 0,1,2,\dots.
\]
The measure is unique because entries in $J_{III}$ are bounded. We rewrite what have been argued as follows. For any $k = 1,2,\dots,$ as $\beta N \to 2c$,
\[
	\bra{L_N, x^k} \to \bra{\nu_{a,b,c}, x^k} \quad \text{almost surely.}
\]
That is to say, in the considering regime, each moment of $L_N$ converges to the corresponding moment of $\nu_{a,b,c}$. This implies the weak convergence of probability measure because the limiting measure is determined by moments (see \cite{Duy-2018} for example). Therefore, we have just proved the following.
\begin{theorem}\label{thm:main}
Let $a, b > -1$ and $c \ge 0$ be fixed. Then in the regime where $\beta N \to 2c$, the sequence of the empirical distribution $L_N$ converges weakly to $\nu_{a, b, c}$, almost surely.
\end{theorem}

\section{Beta Jacobi processes}
Consider the so-called
beta Jacobi processes which are defined to be processes $\{\lambda_i(t)\}_{i = 1, \dots, N}$ in 
\[
 W:=\Big\{ \{x_i\} \in \R^N : 0 \le x_1 \le x_2 \le \cdots \le x_N \le 1 \Big\}
\] 
satisfying the following system of stochastic differential equations (SDEs)
\begin{equation}\label{SDEs}
\begin{cases}
	d \lambda_i =  \sqrt{2\lambda_i (1 - \lambda_i)} db_i + \left(a + 1 - (a+b+2) \lambda_i + \frac\beta2\sum_{j: j \neq i} \frac{2 \lambda_i(1 - \lambda_i)} {\lambda_i - \lambda_j} \right) dt, \\
	\lambda_i(0) = \lambda_0^{(N, i)}, 
\end{cases}
\end{equation}
$i=1, \dots, N$ with initial data $\{\lambda_0^{(N, i)}\} \in W$,
where $a > -1, b > -1,$ and $\beta > 0$. Here $\{b_i\}_{i=1}^N$ are standard Brownian motions. For $\beta = 1$, they are introduced in \cite{Doumerc-2005} as the eigenvalue process of the real Jacobi matrix process. The case for general $\beta > 0$ was then introduced in \cite{Demni-2010} in a relation with radial Dunkl processes. When $\beta \ge 1$, the beta Jacobi processes never collide \cite{Graczyk-Malecki-2014}. They do collide when $0 < \beta < 1$, but the set of $t$ such that $\lambda_i(t) = \lambda_j(t)$, for some $i \neq j$, has Lebesgue measure zero, almost surely. These processes are related to beta Jacobi ensembles in the sense that the following beta Jacobi ensemble restricted in $W$,
\begin{equation}\label{BJE-W}
	const \cdot \prod_{i < j}|\lambda_j - \lambda_i|^\beta \prod_{l = 1}^N \lambda_l^a (1-\lambda_l)^b, \quad 0 \le  \lambda_1 \le \lambda_2 \le \cdots \le \lambda_N \le 1,
\end{equation}
is their stationary distribution. Here is our result on beta Jacobi processes.

We now study the limiting behavior of the empirical measure process $\mu_t^{(N)}$, 
\[
	\mu_t^{(N)} = \frac1N \sum_{i = 1}^N \delta_{\lambda_i(t)},
\]
in the regime where $\beta N \to 2c \in [0, \infty)$. For simplicity, let $c$ be fixed and $\beta = 2c/N$. (We cannot do it when $c = 0$, but the argument is the same.)
For $f \in C^2([0,1])$, by It\^o's formula, we deduce that
\begin{align*}
	d \bra{\mu_t^{(N)}, f} &= \frac1N \sum_{i = 1}^N \sqrt{2 \lambda_i(1 - \lambda_i)} f'(\lambda_i) db_i \\
	&+  \frac1N \sum_i f'(\lambda_i) \left( a + 1 - (a + b + 2) \lambda_i + \frac c N \sum_{j : j \neq i} \frac{2\lambda_i(1 - \lambda_i)}{\lambda_i - \lambda_j} \right) dt \\
	&+ \frac1N \sum_{i = 1}^N \lambda_i(1- \lambda_i) f''(\lambda_i) dt,
\end{align*}
which can be rewritten as 
\begin{align}
	d \bra{\mu_t^{(N)}, f} &= dM_t^{(N)}  +  \bra{\mu_t^{(N)}, (a + 1) f' - (a + b + 2) x f' + x(1-x) f''(x)} dt \notag\\
	&\quad + c \iint \frac{x(1-x) f'(x) - y(1-y) f'(y)}{x - y} d\mu_t^{(N)}(x) d\mu_t^{(N)}(y) dt  \notag\\
	&\quad - \frac cN \bra{\mu_t^{(N)}, [x(1-x)f'(x)]'} dt, \label{Ito-f}
\end{align}
where $M_t^{(N)}$ is a martingale with the quadratic variation 
\[
	[M^{(N)} ]_t = \frac1N\int_0^t \bra{\mu_s^{(N)}, 2 x (1-x) f'(x)^2} ds.
\]
To be more precise, the above identity is true when $\{\lambda_i(t)\}$ are distinct, which holds for almost every $t$.

The moment method has been used in \cite{Trinh-Trinh-2020} to study the limiting behavior of the empirical measure processes of beta Laguerre processes. We observe that the approach is applicable for our model without any difficulty. Therefore, rather than providing detailed proof, we will only sketch some main ideas.
For $k = 1,2, \dots,$ let 
\[	
	S_k^{(N)}(t) = \bra{\mu_t^{(N)}, x^k} = \frac1N \sum_{i = 1}^{N}  \lambda_i(t)^k
\] 
be the $k$th moment process of $\mu_t^{(N)}$. Then equation \eqref{Ito-f} with $f = x^k$ can be expressed in the integral form as follows 
\begin{align}
	S_k^{(N)} (t) &= \bra{\mu_0^{(N)}, x^k} + M_k^{(N)}(t)  - k (2c + a + b + k + 1) \int_0^t S_{k}^{(N)} (s)ds \notag \\
	&\quad +   k (a + k ) \int_0^t S_{k - 1}^{(N)}(s) ds \notag\\
	&\quad + c k \int_0^t \sum_{i = 0}^{k - 1} S_{i}^{(N)}(s) S_{k - 1 - i}^{(N)}(s) ds -  c k \int_0^t \sum_{j = 1}^{k-1} S_{j}^{(N)}(s) S_{k - j}^{(N)}(s) ds \notag\\
	&\quad - \frac cN \int_0^t (k^2 S_{k-1}^{(N)}(s) - k(k+1)S_k^{(N)}(s))ds \label{Ito-sk}.
\end{align}
Here $M_k^{(N)}$ denotes the corresponding martingale part. Given that $\{S_l^{(N)}(t) \}_{0 \le l \le k -1}$ and $M_k^{(N)}$ are known, the equation \eqref{Ito-sk} becomes an initial value ordinary differential equation (ODE) of $\int_0^t S_k^{(N)}(s) ds$, and thus, an explicit formula for $S_k^{(N)}(t)$ can be derived. For fixed $T > 0$, let $\X$ be the space $\cC([0, T], \R)$ of continuous functions on $[0, T]$ endowed with the uniform norm. Then $S_k^{(N)}$ and $M_k^{(N)}$ are random elements on $\X$. By Doob’s martingale inequality, the martingale part $M_k^{(N)}$ is easily shown to converge in probability to zero in $\X$. Then imitate arguments used in \cite{Trinh-Trinh-2020}, we arrive at the following. 
\begin{theorem}\label{thm:dynamical}
Assume that the initial measure $\mu_0^{(N)} = N^{-1} \sum_{i = 1}^N \delta_{\lambda_0^{(N,i)}}$ converges weakly to a probability measure $\mu_0$. Then the following hold.

\medskip
\noindent {\rm(i)}	For each $k = 1, 2, \dots,$ the sequence of $S_k^{(N)}$, as random elements on $\X$, converges in probability to a deterministic limit $m_k(t)$, which is defined inductively as the solution of the following initial value ODE
\begin{align}
	m_k'(t) &= -k (2c + a + b + k + 1) m_k(t) + k (a + k ) m_{k - 1}(t) \notag\\
	&\quad + ck \sum_{i = 0}^{k-1} m_i(t) m_{k - 1 - i}(t) - ck \sum_{j = 1}^{k-1} m_j(t) m_{k  - j}(t), \notag\\
	m_k(0) &=  \lim_{N \to \infty} \bra{\mu_0^{(N)}, x^k}.  \label{ODE-mk}
\end{align}
Here $m_0(t) \equiv 1$.

\medskip 
\noindent {\rm(ii)} Let $\{u_k\}$ be a sequence defined as $u_0 = 1$,
\begin{align}
	u_k = \frac{1}{2c + a + b + k + 1} \bigg( (a + k) u_{k - 1} + c\sum_{i = 0}^{k-1} u_i u_{k - 1 - i} - c \sum_{j = 1}^{k-1} u_j u_{k  - j} \bigg),
\end{align}
for $k = 1,2,\dots$. Then $\lim_{t \to \infty} m_k(t) = u_k$.

\medskip 
\noindent {\rm(iii)} For any $t\ge 0$,  let $\mu_t$ be the unique measure with moments $\{m_k(t)\}$. Then the sequence of the empirical measure process $\mu_t^{(N)}$ converges in probability to $(\mu_t)_{0 \le t \le T}$ as random elements on the space $\cC([0,T], \cP([0,1]))$. Here $\cC([0, T], \cP([0,1]) )$ is the space of continuous maps from $[0,T]$ to the space $\cP([0,1])$ of probability measures on $[0,1]$ endowed with the uniform topology. Moreover, the continuous probability measure-valued process $\mu_t$ converges weakly to $\nu_{a,b,c}$ as $t \to \infty$.

\end{theorem}

	We omit the proof but give some comments. 
\begin{itemize}
\item[(i)] For our model, since all probability measures involved here are supported in $[0,1]$, it follows that the weak convergence of probability measures is equivalent to the convergence of moments. 

\item[(ii)] We have assumed that the initial state $\{\lambda_0^{(N,i)}\} \in W$ is non-random, but the theorem still holds when the initial state is random. In particular, when $\{\lambda_0^{(N,i)}\} \in W$ is distributed as the beta Jacobi ensemble~\eqref{BJE-W}, which is the stationary distribution of the system of SDEs~\eqref{SDEs}, the assumption on the initial measure in the theorem is still fulfilled. This implies that the sequence $\{u_k\}$ coincides with the moments of the limiting measure $\nu_{a,b,c}$. Therefore, the measure $\mu_t$ converges weakly to $\nu_{a,b,c}$ as $t \to \infty$. This proves the last statement in Theorem~\ref{thm:dynamical}.
\end{itemize}

\section{Associated Jacobi polynomials}
In this section, we introduce the classical theory of associated Jacobi polynomials and derive some properties of the limiting measure $\nu_{a, b, c}$. 
Let us recall from the introduction the relation between orthogonal polynomials and Jacobi matrices. Let $\{P_n = x^n + \text{lower powers}\}_{n \ge 0}$ be a sequence of monic polynomials satisfying the following three term recurrence relation
\[
	x P_n = P_{n + 1} + a_{n + 1}P_n + b_n^2 P_{n - 1}, \quad n \ge 0,
\]
where $b_0 = 0, b_n > 0, a_n \in \R, n = 1, 2, \dots$. We form an infinite symmetric tridiagonal matrix $J$, called a Jacobi matrix 
\[
	J = \begin{pmatrix}
		a_1	&b_1	\\
		b_1	&a_2		&b_2\\
		&\ddots	&\ddots	&\ddots
	\end{pmatrix}.
\]
Then there is a probability measure $\mu$ on $\R$ such that 
\begin{equation}\label{spectral-measure}
	 \int x^k d\mu  = J^k(1,1) ,\quad k = 0,1,\dots.
\end{equation}
Moreover, the polynomials $\{P_n\}$ are orthogonal with respect to $\mu$, 
\[
			\int_{\R} P_m(x) P_n(x) d\mu(x) = \begin{cases} 0, &{m \neq n},\\
			b_1^2\cdots b_n^2, &{m = n}.		
			\end{cases}
\]

In case the measure $\mu$ satisfying the moments condition~\eqref{spectral-measure} is unique, or the measure $\mu$ is determined by moments, it is called the spectral measure of the Jacobi matrix $J$. We also call it the probability measure of orthogonal polynomials $\{P_n\}$. The matrix $J$ is referred to as the Jacobi matrix of $\{P_n\}$ as well. When the parameters $\{a_n\}$ and $\{b_n\}$ are bounded, then clearly the measure $\mu$ is unique and has compact support. Assume from now that the measure $\mu$ is unique. Then the polynomials are dense in $L^2(\R, \mu)$, and thus the normalization $p_n = P_n/(b_1\cdots b_n), n \ge 1; p_0 = 1$ which satisfies
	\[
		x p_n = b_{n + 1}p_{n + 1} + a_{n + 1}p_n + b_n p_{n - 1}, \quad n \ge 0, (b_0 :=0)
	\]
becomes an orthonormal basis in $L^2(\R, \mu)$. The matrix $J$ is then the matrix of the linear operator $f \mapsto x f(x)$ on $L^2(\R, \mu)$ with respect to the basis $\{p_n\}$.

Let $S_J(z)$ be the Stieltjes transform of $\mu$,
	\[
		S_J(z) = \int \frac{ d\mu(x)}{x - z} = (J - z)^{-1} (1,1), \quad z \in \C \setminus \R.
	\] 
	It coincides with the $m$-function in the theory of Jacobi matrices.
	Let $J^+$ be the Jacobi matrix obtained from $J$ by removing the first row and the first column. Assume further that the spectral measure of $J^+$ is also unique. Then the following relation holds
\begin{equation}\label{m-function}
		-\frac{1}{S_{J}(z)} = z - a_1 + b_1^2 S_{J^+}(z).
\end{equation}
The proof of all statements above can be found in  \cite{Deift-1999, Simon-book-2011}.

Back to our topic,
Jacobi polynomials are known to be orthogonal polynomials on $[-1,1]$ with respect to the measure $(1 - x)^{\tilde a}(1-x)^{\tilde b} dx$. However, in this work, we consider their variants which are orthogonal on $[0,1]$ with respect to  $x^a (1-x)^bdx$. From the three term recurrence formula for Jacobi polynomials, it turns out that their Jacobi matrix is given by 
\[
	J=\begin{pmatrix}
		\lambda_0 + \mu_0	& \sqrt{\lambda_0 \mu_1}	\\
		 \sqrt{\lambda_0 \mu_1}	&\lambda_1 + \mu_1		& \sqrt{\lambda_1 \mu_2}	\\
		 &\ddots	&\ddots	&\ddots	
	\end{pmatrix}, \quad \text{with $c = 0$}.
\]
Here for convenience, we recall the notations
\begin{align*}
	\hat\lambda_0 &=   \frac{c + a + 1}{2c + a + b + 2},\\
	\lambda_n &=   \frac{n + c + a + 1}{2n + 2c + a + b + 2} \frac{n + c + a + b + 1}{2n + 2c + a + b + 1}, \quad n \ge 0,\\
	\mu_n &=   \frac{n + c}{2n + 2c + a + b + 1} \frac{n + c + b}{2n + 2c + a + b}, \quad n \ge 0.
\end{align*}

\emph{Associated Jacobi polynomials: Model I.} Orthogonal polynomials related to the following (Jacobi) matrix was considered in  \cite{Wimp-1987}
\[
	J_{I}=\begin{pmatrix}
		\lambda_0 + \mu_0	& \sqrt{\lambda_0 \mu_1}	\\
		 \sqrt{\lambda_0 \mu_1}	&\lambda_1 + \mu_1		& \sqrt{\lambda_1 \mu_2}	\\
		 &\ddots	&\ddots	&\ddots	
	\end{pmatrix}.
\]
Here the requirement is that $c \ge 0, c+a > 0, c+b>0$. Note that for $c \in \N$, $\{\lambda_n(c), \mu_n(c)\} = \{\lambda_{n + c}(0), \mu_{n + c}(0)\}$ are shifted from the parameters of Jacobi polynomials, and hence the name. Explicit representations for orthogonal polynomials and for the spectral measure were derived.

\emph{Associated Jacobi polynomials: Model II.} With a slight modification of the first entry, another model of associated Jacobi polynomials was studied in \cite{Ismail-Masson-1991}. In the terminology here, this model deals with the following matrix 
\[
	J_{II}=\begin{pmatrix}
		\lambda_0 	& \sqrt{\lambda_0 \mu_1}	\\
		 \sqrt{\lambda_0 \mu_1}	&\lambda_1 + \mu_1		& \sqrt{\lambda_1 \mu_2}	\\
		 &\ddots	&\ddots	&\ddots	
	\end{pmatrix},
\]
which was motivated from the study of birth and death processes in which $\mu_0$ representing the death rate at $0$ is naturally assumed to be zero.

\emph{Associated Jacobi polynomials: Model III.} In this paper, we encounter with a new type of associated Jacobi polynomials, which is different from Model II by replacing $\lambda_0$ with $\hat \lambda_0$
\[
	J_{III}=\begin{pmatrix}
		\hat\lambda_0 	& \sqrt{\hat\lambda_0 \mu_1}	\\
		 \sqrt{\hat\lambda_0 \mu_1}	&\lambda_1 + \mu_1		& \sqrt{\lambda_1 \mu_2}	\\
		 &\ddots	&\ddots	&\ddots	
	\end{pmatrix}.
\]
Here recall that $c \ge 0$ and $a, b  > -1$. However, for the meaning of this model, we only need $c + 1 > 0, c + a + 1>0, c + b + 1 > 0$. Observe that when the first row and the first column of $J_{III}$ are removed, we get Model I (with parameter $c + 1$). Thus, by using the relation~\eqref{m-function}, we see that
\[
	S_{III}(z) = -\frac{1}{z - \hat \lambda_0 + \hat \lambda_0 \mu_1 S_{I}(z; c+1)},
\]
where $S_{III}(z)$ is the Stieltjes transform of the spectral measure of $J_{III}$, and $S_{I}(z; c)$ denotes the Stieltjes transform of the spectral measure in Model I. It was shown in \cite{Wimp-1987} that 
\[
	S_{I}(z; c) = - \frac{_2F_1(c+1, c+a+1; 2c + a + b + 2; 1/z)}{z \,_2F_1(c, c+a; 2c + a + b; 1/z)}.
\]
Here $_2F_1$ denotes the hypergeometric function.
Then after some simplifications, we deduce that 
\[
	S_{III}(z) = - \frac{_2F_1(c+1, c+a+1; 2c + a + b + 2; 1/z)}{z \,_2F_1(c, c+a+1; 2c + a + b +2; 1/z)}.
\]
From which, an explicit formula for $\nu_{III}$ can be derived.

For the sake of the completeness, we include here the formulas of the Stieltjes transforms  in all three models.
\begin{theorem}
Let $S_{i}, i \in \{I, II, III\}$ be the Stieltjes transforms of the probability measures of associated Jacobi polynomials in Model $i$. Then the following hold
	\begin{align*}
		S_{I}(z) = - \frac{_2F_1(c+1, c+a+1; 2c + a + b + 2; 1/z)}{z \,_2F_1(c, c+a; 2c + a + b; 1/z)};\\
		S_{II}(z) = - \frac{_2F_1(c+1, c+a+1; 2c + a + b + 2; 1/z)}{z \,_2F_1(c, c+a+1; 2c + a + b +1; 1/z)};\\
		S_{III}(z) = - \frac{_2F_1(c+1, c+a+1; 2c + a + b + 2; 1/z)}{z \,_2F_1(c, c+a+1; 2c + a + b +2; 1/z)}.
	\end{align*}
\end{theorem}

Next, we imitate the method in \cite{Ismail-Masson-1991} to derive an explicit density for the spectral measure in Model III. This approach requires $a$ is not an integer number.

{\bf Step 1.} 
Let $\{R_n = R_n^{a,b}(x; c)\}$ be defined as
\[
\frac{n +  c + 1}{n + c + a + 1}\lambda_n R_{n + 1} = (x - \lambda_n - \mu_n) R_n - \frac{n + c + a}{n + c}\mu_n R_{n - 1}, \quad n \ge 0,
\]
with initial conditions
\[
	R_{-1}=0, \quad R_0 = 1.
\]
Note that $\{\hat R_n := R_{n-1}^{a,b}(x; c+1)\}$ satisfies the same recurrence relation with $\{R_n\}$ but with different initial conditions
\[
	\hat R_0 = 0, \quad \hat R_1 = 1.
\]
The following formula for $R_n$ was known \cite{Wimp-1987} (Eq.~(28))
\begin{align*}
	&R_n=R_n^{(a,b)}(x;c)\\
	&=\frac{(-1)^n \Gamma(c + 1) \Gamma(\gamma + c)}{a \Gamma(a + c) (\gamma + 2c - 1) \Gamma(\gamma + c - a - 1)}\\
	&\times \bigg\{ \frac{\Gamma(\gamma + c - a - 1) \Gamma(n + a + c + 1)}{\Gamma(\gamma + c - 1)\Gamma(n + c + 1) } \\
	&\qquad\qquad \times \,_2F_1(c, 2 - \gamma -c; 1-a; x) \,_2F_1(-n-c, n+\gamma+c;a + 1; x)\\
	&\qquad - \frac{\Gamma(a + c) \Gamma(n + \gamma + c - a)}{\Gamma(c) \Gamma(n + c + \gamma)} \\
	&\qquad\qquad \times \,_2F_1(1-c, \gamma + c - 1; a + 1; x) \,_2F_1(n + c + 1, 1-n-\gamma-c; 1-a; x) \bigg\},
\end{align*}
with $\gamma = a + b + 1$.

We consider the following sequence of orthogonal polynomials $\{P_n\}$ 
\begin{align*} 
	&P_0 = 1,\quad 
	\frac{c + 1}{c + a + 1} \hat\lambda_0 P_{1} = (x - \hat\lambda_0 ) P_0,\\
	&\frac{n +  c + 1}{n + c + a + 1}\lambda_n P_{n + 1} = (x - \lambda_n - \mu_n) P_n - \frac{n + c + a}{n + c}\mu_n P_{n - 1}, \quad n \ge 1.
\end{align*}
Since $\{P_n\}$, $\{R_n\}$ and $\{\hat R_n\}$  satisfy the same recurrence relation and since $\{R_n\}$ and $\{\hat R_n\}$ are linearly independent, the sequence $\{P_n\}$ can be expressed as 
\[
	P_n = A R_n + B \hat R_n,
\]
where $A$ and $B$ are constants not depending on $n$, which can be solved from the initial conditions. As the result, we get the following expression for $\{P_n\}$
\[
	P_n = R_n^{a,b}(x; c) + \bigg\{\frac{c(c+b)(2c + \gamma + 1)}{(c+1)(2c + \gamma - 1)(c+ \gamma)} - \frac{c(2c + \gamma + 1)}{(c+1)(c + \gamma)} x  \bigg\} R_{n - 1}^{a,b}(x; c+1).
\]

\begin{lemma}
The polynomial $P_n$ has the following expression
\begin{align*}
	(-1)^n P_n &=  \frac{\Gamma(c+1) \Gamma(n+c+a+1)}{\Gamma(n+c+1)\Gamma(c+a+1)}\,_2F_1(c,-c-\gamma;-a;x) \\
	&\qquad \qquad\times  \,_2F_1(-c-n,c+n+\gamma;1+a;x) \\
	&\quad -\frac{c\Gamma(\gamma + c + 1) \Gamma(n+c+b+1)}{a(a+1)\Gamma(\gamma + n+c)\Gamma(c+b+1)} x(1-x) \,_2F_1(1-c,1+c+\gamma;2+a;x)\\
	&\qquad \qquad\times \,_2F_1(1+c+n, -c-n-a-b;1-a;x).
\end{align*}
\end{lemma}

{\bf Step 2.} It follows from the definition of $\{P_n\}$ that $P_n$ is a polynomial of degree $n$ with the highest coefficient 
\[
	\bigg(\hat \lambda_0 \lambda_1 \cdots \lambda_{n - 1} \frac{(c+1)_n}{(c+a+1)_n}  \bigg)^{-1}.
\]
Here $(q)_n$ is the Pochhammer symbol defined by
\[
	(q)_n = \begin{cases}
		1, & n = 0, \\
		q(q+1)\cdots (q + n - 1), & n \ge 1.
	\end{cases}
\]
Let $\nu=\nu_{a,b,c}$ be the spectral measure of $J_{III}$. Then the sequence $\{p_n = P_n / \zeta_n\}$ becomes an orthonormal system in $L^2(\R, \nu)$, where 
\[
	\zeta_n = \bigg(\frac{\mu_1 \mu_2 \cdots \mu_n}{\hat \lambda_0 \lambda_1 \cdots \lambda_{n - 1}}\bigg)^{1/2} \frac{(c+a+1)_n}{(c+1)_n}.
\]
By using the following asymptotic of the gamma function
\[
	\frac{\Gamma(a+n)}{\Gamma(b + n)} \approx n^{a - b}\quad \text{as } n \to \infty,
\]
we obtain the asymptotic of $\zeta_n$ as 
\[
	\zeta_n \approx \frac{1}{\sqrt{2n}} \bigg(\frac{\Gamma(c+1) \Gamma(c+ a+b + 2)}{\Gamma(c+a+1)\Gamma(c+b+1)} \bigg)^{1/2}.
\]
Here $f(n) \approx g(n)$ as $n \to \infty$ means that $f(n) / g(n) \to 1$ as $n \to \infty$.
In addition, by using the asymptotic of the hypergeometric function $\,_2F_1$
\begin{align*}
	&\, _2F_1(a+n, b - n; c; \sin^2\theta) \\
	&\quad \approx \frac{\Gamma(c) n^{-c + 1/2} (\cos \theta)^{c-a-b-1/2}}{\sqrt{\pi} (\sin \theta)^{c-1/2}} \cos \Big[2 n \theta + (a-b)\theta + \frac{\pi}{2}\Big(\frac 12 - c \Big) \Big]
\end{align*}
as $n \to \infty, \theta \in (0, \pi/2)$, we can derive the asymptotic of $p_n(x)$ for $x \in (0,1)$ as
\[
	\limsup_{n \to \infty} p_n^2(x)  = \frac{2}{\pi} \frac{\Gamma(c+a+1)\Gamma(c+b+1)}{\Gamma(c+1) \Gamma(c+ a+b + 2)} (1-x)^{-b-1/2} x^{-a-1/2} |U(x) + e^{i \pi a} V(x)|^2,
\]
where 
\begin{align}
	U(x)&=\frac{\Gamma(c+1) \Gamma(a + 1)}{\Gamma(1+c+a)}  \,_2F_1(c, -c-a-b - 1; -a; x), \label{U}\\
	V(x)&=\frac{-\pi c \Gamma(c+a+b+2)}{\sin(\pi a) \Gamma(1+c+b)\Gamma(2+a)}(1-x)^{1+b} x^{1+a}  \,_2F_1(1-c, 2+c+a+b, 2+a, x).\label{V}
\end{align}

{\bf Step 3.} We now use the following result to derive the density of $\nu$.
\begin{lemma}[{\cite[pp.\ 141--143]{Nevai-1979}}]
Let $J$ be a Jacobi matrix
\[
	J = \begin{pmatrix}
		a_1		&b_1	\\
		b_1		&a_2		&b_2	\\
		&\ddots 	&\ddots 	&\ddots
	\end{pmatrix},
	\quad a_n \in \R, b_n > 0
\]
with bounded parameters and let $\nu$ be its spectral measure. Let $\{p_n(x)\}$ be its orthonormal polynomials. Assume that 
\[
	\sum_{n = 1}^\infty	(|a_n - a| + |b_n - b|) < \infty,
\]
where $a \in \R$ and $b > 0$ are constants. Then for almost every $x \in (a-2b, a+2b)$,
\[
	\limsup_{n \to \infty} \sqrt{(x-a+2b)(a+2b-x)}p_n^2(x)\nu_{ac}(x) = \frac{2}{\pi}. 
\]
Here $\nu_{ac}(x)$ is the density of the absolutely continuous part of $\nu$. Moreover, $\supp(\nu_{ac}) = [a-2b, a+2b]$ and the singular part is supported outside $(a-2b, a+2b)$.
\end{lemma}

Using the lemma, we conclude that the density of the absolutely continuous part of the spectral measure $\nu$ is given by
\[
	\nu_{ac}(dx) = \frac{\Gamma(c+1) \Gamma(c+ a+b + 2)}{\Gamma(c+a+1)\Gamma(c+b+1)} \frac{ x^{a}(1-x)^b }{|U(x) + e^{i \pi a} V(x)|^2},\quad 0<x<1.
\]
The singular part is actually zero, which can be proved in a similar way as in Model II \cite{Ismail-Masson-1991}. In conclusion, we obtain the following explicit formula for the density of the spectral measure $\nu_{a,b,c}$ in Model III.
\begin{theorem}
	Assume that $c \ge 0,  c+a > 0$ and $c+b>0$ and $a$ is not an integer number. Then  the spectral measure $\nu = \nu_{a,b,c}$  in Model III is absolutely continuous with density  
	\[
		\nu(dx) = \frac{\Gamma(c+1) \Gamma(c+ a+b + 2)}{\Gamma(c+a+1)\Gamma(c+b+1)} \frac{ x^{a}(1-x)^b }{|U(x) + e^{i \pi a} V(x)|^2},\quad 0<x<1.
	\]
with $U$ and $V$ given in \eqref{U} and \eqref{V}.
\end{theorem}

\bigskip
\noindent{\textbf{Acknowledgements.}} 
This work is supported by JSPS KAKENHI Grant Number JP19K14547 (K.D.T.). The authors would like to thank a referee for helpful comments.

\begin{footnotesize}

\end{footnotesize}


\end{document}